\documentclass[12pt]{article}
\usepackage{amsmath}
\usepackage{amssymb}

\usepackage[cp1251]{inputenc}
\usepackage[russian]{babel}
\newcommand{\il}[2]{\int\limits_{#1}^{#2}}
\newcommand{\ilp}[1]{\int\limits_{#1}^{+\infty}}

\newcommand{\ph}{\phantom{a}}
\newcommand{\phh}{\phantom{aaa}}

\newcommand{\sist}[2]{\left\{
\begin{array}{l}
{#1}\\
\ph\\
{#2}
\end{array}
\right.}

\begin{document}
MSC  34L30,   34C99

\vskip 20pt
{\bf \centerline {Properties of solutions of  matrix}

\centerline{Riccati equations}

\vskip 20pt \centerline {G. A. Grigorian}}
\centerline{\it Institute  of Mathematics NAS of Armenia}
\centerline{\it E -mail: mathphys2@instmath.sci.am}
\vskip 20 pt

Abstract. In this paper we study properties of regular solutions of matrix Riccati equations. 
 The obtained results are used to study the asymptotic behavior of solutions of linear systems of ordinary differential equations.


\vskip 30pt

{\bf 1. Introduction}. Let $P(t), \ph Q(t), \ph R(t)$ and $S(t)$ be complex-valued continuous matrix functions of dimension $n\times n$ on $[t_0,+\infty)$. Consider the matrix Riccati equation
$$
Z' + Z P(t) Z + Q(t) Z + Z R(t) + S(t) = 0, \phh t \ge t_0. \eqno (1.1)
$$
It follows from the general theory of normal systems of ordinary differential equations for any matrix $M$ of dimension $n\times n$ and $t_1 \ge t_0$ there exists $t_2 > t_1 \ph (t_2 \le +\infty)$ such that the solution $Z(t)$ of the last equation exists on $[t_1,t_2)$ and is the unique. Main interest from the point of view of qualitative theory of differential equations represents  the case $t_2 = +\infty$.

{\bf Definition 1.1.} {\it A solution $Z(t)$ of Eq. (1.1) is called $t_1$-regular if it exists on $[t_1,+\infty).$}

Obviously every $t_0$-regular solution of Eq. (1.1) is also $t_1$-regular for every $t_1 > t_0$. For this every $t_0$-regular solution of Eq. (1.1) we will call just regular solution of Eq. (1.1)

{\bf Remark 1.1.} {\it Criteria for existence of regular solutions of Eq. (1.1) are obtained in [1-3] (see [1], Theorems 3.5, 3.6, [2,3]).}

Some interest represents the investigation  properties of $t_1$-regular solutions of Eq. (1.1) in the context of their application in the study of properties  of solutions of linear systems of ordinary differential equations. In the scalar ($n=1$) and the quaternionic cases some properties of $t_1$-regular solutions of Eq. (1.1) are studied in [4-9]  and are used to qualitative study of some linear ordinary differential equations and systems of such equations (see [10-16]).

In this paper we study some properties of $t_1$-regular  solutions of Eq. (1.1). The obtained results we use to study the asymptotic behavior of solutions of linear systems of ordinary differential equations.
\vskip 20pt

{\bf 2. Auxiliary relations.} Let $Z(t)$ be a solution of Eq. (1.1) on $[t_1,t_2) \ph (t_0\le t_1 < t_2 \le +\infty)$. Then any solution $Z_1(t)$ of Eq. (1.1) on $[t_1,t_2)$ is given by the formula (see [17], pp. 139, 140, 158, 159, Theorem 6.2)
$$
Z_1(t) = Z(t) + \{\phi_Z(t)\Lambda^{-1}(t_1)[I + \Lambda(t_1)\mu_Z(t_1,t)]\psi_Z(t)\}^{-1}, \phh t \in [t_1,t_2), \eqno (2.1)
$$
where $\Lambda(t_1) \equiv Z_1(t_1) - Z(t_1), \ph (det \hskip 2pt \Lambda(t_1) \ne 0)\ph \phi_Z(t)$ and $\psi_Z(t)$ are solutions of the linear matrix equations
$$
\phi' = [P(t)Z(t) + R(t)]\phi = 0, \phh t \in [t_1,t_2),
$$
$$
\psi' = \psi[Z(t) P(t) + Q(t)] = 0, \phh t \in [t_1,t_2)
$$
respectively with $\phi_Z(t_1) = \psi_Z(t_1) = I, \ph I$ is the identity matrix of dimension $n\times n$,
$$
\mu_Z(t_1,t) \equiv \il{t_1}{t}\phi_Z^{-1}(\tau) P(\tau) \psi_Z^{-1}(\tau) d \tau, \phh t \in [t_1,t_1)
$$
(since $\phi_Z(t_1) = \psi_Z(t_1) = I$ by the Liuoville's formula $det \hskip 2pt \phi_Z(t) \ne 0, \ph det \hskip 2pt \psi_Z(t) \ne~ 0, \linebreak t\in~ [t_1,t_2)$).
Obviously this formula we can rewrite in the form
$$
Z_1(t) = Z(t) + \psi_z^{-1}(t)[I + \Lambda(t_1)\mu_Z(t_1,t)]^{-1}\Lambda(t_1)\phi_Z^{-1}(t), \phh t \in [t_1,t_2). \eqno (2.2)
$$
Note that unlike of the formula (2.1) the right hand part of (2.2) is meaning full not only for $\det \Lambda(t_1) \ne 0$ but also for $det \hskip 2pt \Lambda(t_1) = 0$. Show that (2.2) is valid also for $\det  \Lambda(t_1) = 0$. Suppose $\det  \Lambda (t_1) = 0$. Chose $\varepsilon_0 > 0$ so small that for every $\varepsilon \in (0,\varepsilon_0)$ the relation $\det  (\varepsilon I + \Lambda(t_1))\ne 0$ holds (this is possible always, because $\det  (\varepsilon I + \Lambda(t_1))$ is a polynomial with respect to $\varepsilon \in\mathbb{C}$. Let for each $\varepsilon \in (0,\varepsilon_0)$ the matrix function $Z_\varepsilon(t)$ be the solution of Eq. (1.1) with $Z_\varepsilon(t_1) = \varepsilon I + Z_1(t_1)$. Set $\Lambda_\varepsilon(t_1)\equiv Z_\varepsilon(t_1) - Z(t_1)$. Since $\det  \Lambda_\varepsilon(t_1)\ne 0$ according to (2.2) we have
$$
Z_\varepsilon(t) = Z(t) + \psi_z^{-1}(t)[I + \Lambda_\varepsilon(t_1)\mu_Z(t_1,t)]^{-1}\Lambda_\varepsilon(t_1)\phi_Z^{-1}(t) = \phantom{aaaaaaaaaaaaaaaaaaaaaaaa}
$$
$$
\phantom{aaaaaaaaaa}=Z(t) + \psi_z^{-1}(t)[I + (\varepsilon I +\Lambda(t_1))\mu_Z(t_1,t)]^{-1}\Lambda_\varepsilon(t_1)\phi_Z^{-1}(t), \phh t \in [t_1,t_3), \eqno (2.3)
$$
for some $t_3\in (t_1,t_2).$ Let $\delta > 0$ be fixed. By the continuous dependence of solutions of Eq. (1.1) on their initial values we can take $\varepsilon_0$ so small that $Z_\varepsilon(t)$ exists on $[t_1, t_2-\delta]$ for all $\varepsilon \in (0,\varepsilon_0)$ (therefore we can take that $t_3 \ge t_2 - \delta$ and
$$
||Z_\varepsilon(t_1) - Z_1(t_1)|| \le \delta, \phh t \in [t_1,t_2 - \delta]
$$
for all $\varepsilon \in (0,\varepsilon_0)$, where $||\cdot ||$ denotes any matrix norm. This together with (2.3) implies
$$
\Delta(t) \equiv ||Z_1(t) - Z(t) - \psi_z^{-1}(t)[I + +\Lambda(t_1)\mu_Z(t_1,t)]^{-1}\Lambda(t_1)\phi_Z^{-1}(t)|| \le \phantom{aaaaaaaaaaaaaaaaaaaaaa}
$$
$$
\le|Z_1(t) - Z_\varepsilon(t)|| + ||Z\varepsilon(t) - \psi_z^{-1}(t)[I + +\Lambda(t_1)\mu_Z(t_1,t)]^{-1}\Lambda(t_1)\phi_Z^{-1}(t)|| \le \phantom{aaaa}
$$
$$
\phantom{aa}\le \delta + ||\psi_Z^{-1}(t)|| \cdot ||[I + \Lambda_\varepsilon(t_1)\mu_Z(t_1,t)]^{-1}(\Lambda_\varepsilon(t_1) - \Lambda(t_1))|| \cdot||\phi_Z^{-1}(t)|| +
$$
$$
+ ||\psi_Z^{-1}(t)|| \cdot || \{[I + \Lambda_\varepsilon(t_1)\mu_Z(t_1,t)]^{-1} - [I + \Lambda (t_1)\mu_Z(t_1,t)]^{-1}\}||\cdot || \Lambda(t_1)\psi_z^{-1}(t)||,  \eqno (2.4)
$$
$t \in [t_1,t_2 -\delta].$ Obviously
$$
\lim\limits_{\varepsilon \to 0}||\Lambda\varepsilon(t_1) - \Lambda(t_1)|| = 0,
$$
and
$$
\lim\limits_{\varepsilon \to 0}\max\limits_{t \in [t_1,t_2-\delta]}|| \{[I + \Lambda_\varepsilon(t_1)\mu_Z(t_1,t)]^{-1} - [I + \Lambda (t_1)\mu_Z(t_1,t)]^{-1}\}|| = 0,
$$
which together with (2.4) implies $\Delta(t) \le \delta, \ph t \in [t_1,t_2-\delta]$. From here it follows (2.2).

Let $Z_1(t)$ and $Z_2(t)$ be solutions of Eq. (1.1) on $[t_1,t_2) \ph (t_0 \le t_1 < t_2 \le +\infty).$ Set $\Lambda_{jk}(t_1)\equiv Z_J(t_1) - Z_k(t_1), \ph j, \ph k = 1,2.$ By (2.2) we have
$$
\phi_{Z_1}^{-1}(t)P(t)[Z_j(t) - Z_k(t)]\phi_{Z_j}(t) = \phi_{Z_1}^{-1}(t)P(t)\psi_{Z_j}^{-1}(t)[I + \Lambda_{jk}(t)\mu_{Z_j}(t_1,t)]^{-1}\Lambda_{jk}(t_1), \eqno (2.5)
$$
$t\in [t_1,t_2).$ Assume $\det \Lambda_{jk}(t_1) \ne 0$. Then from (2.5) we obtain
$$
\Lambda_{jk}(t_1)\phi_{Z_1}^{-1}(t)P(t)[Z_j(t) - Z_k(t)]\phi_{Z_j}(t)\Lambda^{-1}_{jk}(t_1) = \phantom{aaaaaaaaaaaaaaaaaaaaaaaaaaaaaaaaaaaaaaa}
$$
$$
\phantom{aaaaaaaaaaaaaaaaaa}=\Lambda_{jk}(t_1)\phi_{Z_1}^{-1}(t)P(t)\psi_{Z_j}^{-1}(t)[I + \Lambda_{jk}(t_1)\mu_{Z_j}(t_1,t)]^{-1}, \phh t \in [t_1,t_2).
$$
Hence,
$$
[I + \Lambda_{jk}(t_1)\mu_{Z_j}(t_1,t)]' = A_{Z_j,Z_k}(t)[I + \Lambda_{jk}(t_1)\mu_{Z_j}(t_1,t)] \phh t \in [t_1,t_2),
$$
where
$$
A_{Z_j,Z_k}(t) \equiv \Lambda_{jk}(t_1)\phi_{Z_1}^{-1}(t)P(t)[Z_j(t) - Z_k(t)]\phi_{Z_j}(t)\Lambda^{-1}_{jk}(t_1), \phh t \in [t_1,t_2).
$$
By virtue of the Liouville's formula from here we obtain
$$
\det [I + \Lambda_{jk}(t_1)\mu_{Z_j}(t_1,t)] = \exp\biggl\{\il{t_1}{t} tr \hskip 2pt \Bigl[P(\tau)(Z_j(\tau) - Z_k(\tau))\Bigr]d \tau\biggr\}, \eqno (2.6)
$$
$\in [t_1,t_2), \ph j,k =1,2.$
It is not difficult to verify that the obtained equality is valid also for $\det \Lambda_{j,k}(t_1) \ne 0$
(to prove this statement it is enough to replace  $Z_j(t)$ in (2.6) by $\varepsilon I + Z_j(t), \ph \varepsilon \in (0,\varepsilon_0)$, where $\varepsilon > 0$ is  such  small that $\det [\varepsilon I + Z_j(t_1)] \ne 0, \ph \varepsilon \in~ (0,\varepsilon_0)$ and using estimates, similar to (2.4), take limit, when $\varepsilon \to 0$). From (2.6) we obtain immediately
$$
\det \{[I + \Lambda_{jk}(t_1)\mu_{Z_j}(t_1,t)] \hskip 1pt [I + \Lambda_{kj}(t_1)\mu_{Z_k}(t_1,t)]\} \equiv 1, \phh t\in [t_1,t_2). \eqno (2.7)
$$
\vskip 20pt

{\bf 3. Properties of solutions of Eq. (1.1).} Let $t_1 \ge t_0.$

{\bf Definition 3.1.} {\it A $t_1$-regular solution $Z(t)$ of Eq. (1.1) is called $t_1$-normal if there exists a neighborhood $U(Z(t_1))$ of $Z(t_1)$ such that every solution $\widetilde{Z}(t)$ of Eq. (1.2) with $\widetilde{Z}(t_1) \in U(Z(t_1))$ is also $t_1$-regular, otherwise $Z(t)$ is called $t_1$-extremal.}

{\bf Definition 3.2.} {\it Eq. (1.1) is called regular if it has at least one regular solution.}

{\bf Remark 3.1.} {\it Since the solutions of Eq. (1.1) are continuously dependent on their initial values every $t_1$-normal ($t_1$-extremal) solution of Eq. (1.1) is also a $t_2$-normal ($t_2$-extremal) solution of Eq. (1.1) for all $t_2 > t_1$. Due to this a $t_1$-normal ($t_1$-extremal) solution of Eq. (1.2) we will just call a normal (a extremal) solution of Eq. (1.1). Note that a $t_2$-normal ($t_2$-extremal) solution of Eq. (1.1) may not be a $t_1$-normal ($t_1$-extremal) solution of Eq.(1.1) if $t_1 < t_2$, because a $t_2$-regular solution of Eq. (1.1) may not be $t_1$-regular for $t_1< t_2$.}

Denote by $Q(t,t_1,\Lambda)$ the general solution of Eq. (1.1) in the region $G_{t_1} \equiv \{(t,Z) : \ph t \in~ I_{t_1}(\Lambda), \ph  Z, \Lambda \in\mathbb{C}^{n^2}\},$ where $I_{t_1}(\Lambda)$ is the maximum existence interval for the solution $Z(t)$ of Eq. (1.1) with $Z(t_1) = \Lambda$.

{\bf Example 3.1.} {\it Consider the equation
$$
Z' + Z P(t) Z = 0, \phh t \ge -1. \eqno (3.1)
$$
According to (2.2) the general solution of this equation in the region $G_0 \cap [-1, +\infty) \times \mathbb{C}^{n^2}$ is given by formula
$$
Q(t,0,\lambda) = \Bigl[I + \Lambda\il{0}{t}P(\tau)  d \tau\Bigr]^{-1}\Lambda, \ph \Lambda \in \mathbb{C}^{n^2}, \ph \det \Bigl[I + \Lambda\il{0}{t}P(\tau)  d \tau\Bigr] \ne 0,  t \in I_0(\Lambda). \eqno (3.2)
$$
Assume $P(t)$ has a bounded support. Then from (3.2) is seen that Eq. (3.1) has no $0$-extremal solution, and all its solutions $Q(t,0,\lambda)$ with  $||\lambda||< \sup\limits_{t \ge 0}||\il{0}{t}P(\tau) d \tau||$ are\linebreak $0$-~normal. If $P(t) = p(t) I, \ph t \ge 0$, where $p(t)$ is a non negative function with an unbounded support and $J_0 \equiv \ilp{0}p(\tau) d \tau < +\infty$ then from (3.2) is seen that the solution $Z_*(t) = Q(t,0,-\frac{1}{J_0}I)$ is $0$-extremal, all the solutions $Z_\lambda(t) \equiv Q(t,0,\lambda I), \ph t \ge 0, \ph \lambda \in \mathbb{R})$ are $0$-normal, if $\lambda < -\frac{1}{J_0}$, and are not $0$-regular if $\lambda > -\frac{1}{J_0}$.

Let $(U_{ij}(x_{11},\dots,x_{1n},x_{21},\dots,x_{2n},\dots,x_{n1},\dots,x_{nn}))_{i,j=1}^n$ be a continuously differentiable mapping from $[-1,1]^{n^2}$ onto the unite sphere $\mathbb{S}^{n^2}$ in $\mathbb{R}^{n^2}\subset\mathbb{C}^{n^2}$, and let $\lambda_{ij} \in \mathbb{R}, \ph i,j=\overline{1,n}$ be any rationally independent numbers. Assume
$$
\il{0}{t}P(\tau) d \tau = \frac{2}{\pi}\arctan t (U_{ij}(\sin \lambda_{11} t, \dots, \sin \lambda_{nn} t))_{i,j=1}^n, \phh t \ge t_0.
$$
Then from (3.2) is seen that the solutions $Q(t,0,\Lambda)$ of Eq. (3.1) with $||\Lambda|| <1$ are $0$-normal, and since the set
$$
\{ (U_{ij}(\sin \lambda_{11} t, \dots, \sin \lambda_{nn} t))_{i,j=1}^n, \ph t \ge 0\}
$$
is everywhere dense in $\mathbb{S}^{n^2}$ (see [18]) the solutions of Eq. (3.1) with $\Lambda \in \mathbb{R}^{n^2}, \ph ||\Lambda|| =1$ are $0$-extremal.

{\bf Example 3.2.} {\it Let $u_0 \in \mathbb{C}^{n^2}, \ph \det u_0 \ne 0$ and let $0 < r < R < +\infty$. Denote $K_{r,R}(u_0) \equiv \{Z \in \mathbb{C}^{n^2} : r < |Z - u_0| < R\}$ - an annulus in $\mathbb{C}^{n^2}$ with a center $u_0$ and radiuses $r$ and $R$. For any $\varepsilon > 0$ denote $K_{\varepsilon, r,R}(u_0) \equiv \{\xi_1, ..., \xi_m \in K_{r,R}(u_0): \ph if \ph u \in K_{r,R}(u_0) \ph then \ph  there \ph exists \ph s \in~ \{1, ..., m\} \ph such \ph that \ph |u - \xi_s| < \varepsilon\}$- a finite $\varepsilon$-net for $K_{r,R}(u_0)$ (here $m$ depends on $\varepsilon$). Consider the sequence of $\frac{1}{2n}$-nets: $\{K_{\frac{1}{2n}, \frac{1}{n}, n}(u_0)\}_{n=1}^{+\infty}$. Let the function $f(t) \equiv \il{0}{t} P(\tau) d \tau, \ph t \ge 0$ has the following properties: $f(t) \ne u_0, \ph t\in [0,1]$; when t varies from $n$ to $n+1 \ph (n=1, 2, ...)$ the curve $f(t)$ crosses all points of $K_{\frac{1}{2n},\frac{1}{n}, n}(u_0)$ (i. e. for every $v \in K_{\frac{1}{2n},\frac{1}{n}, n}(u_0)$ there exists $\zeta_v \in [n,n+1]$ such that $f(\zeta_v) = v); \ph  f(t) \in K_{\frac{1}{2n}, +\infty}(u_0) \ph n =1, 2, .... , \ph t \ge 1.$ From these properties it follows that for every $T\ge 0$ the set $\{f(t) : t\ge T\}$ is everywhere dense in $\mathbb{C}^{n^2}$ and $f(t) \ne u_0, \ph t \ge 0$. Hence from (3.2) it follows that Eq. (3.1) has no $t_1$-normal solutions for all $t_1 \ge 0$ and has at least two extremal solutions: $q_1(t) \equiv 0$ and $q_2(t)$ with $q_2(0) = - u_0^{-1}.$ By analogy using $\frac{1}{2n}$-nets $K_{\frac{1}{2n},\frac{1}{n},n}(u_0;...u_l)\equiv \{\xi_1, ..., \xi_m \in \bigcap\limits_{k=0}^l K_{\frac{1}{n},n}(u_k) : u \in~  \bigcap\limits_{k=0}^l K_{\frac{1}{n},n}(u_k)\Rightarrow \exists s \in \{1, ..., m\} : |u - \xi_s| < \frac{1}{2n}\}$ of the intersections $\bigcap\limits_{k=1}^l K_{\frac{1}{n},n}(u_k)$ in place of $K_{\frac{1}{2n},\frac{1}{n},n}(u_0), \ph n =1, 2, ....$ one can show that there exists a Riccati equation, which has no $t_1$-normal solutions and has at least $l+2$ $t_1$-extremal solutions for all $t_1 \ge 0$.}

Hereafter we will assume that Eq. (1.1) has at least two $t_1$-regular solutions for some $t_1\ge t_0$.

{\bf Theorem 3.1.} {\it A $t_1$-regular solution  $Z(t)$  of  Eq. (1.1) is $t_1$-normal if and only if $\mu_Z(t_1,t), \ph t \ge t_1$ is bounded by $t$ on $[t_1,+\infty)$.}

Proof. Necessity.  Suppose $\mu_Z(t_1,t)=(\mu_{Z,i,j}(t_1,t))_{i,j=1}^n$ is not bounded by $t$ on $[t_1,+\infty)$. Then for some $i=i_0$ and $j=j_0$ there exists an infinitely large sequence $t_1 < t_2 <\dots < t_m < \dots$ such that
$$
\lim\limits_{m\to +\infty}|\mu_{Z,i_0,j_0}(t_1,t_m)| = + \infty. \eqno (3.3)
$$
Let $Z_m(t), \ph m=1,2,\dots$ be solutions of Eq. (1.1) with $Z_m(t_1) = Z(t_1) + A_m, \ph m=1,2,\dots$, where $A_m \equiv (a_{i,j,m})_{i,j=1}^n  \ph m=1,2,\dots$ satisfy the conditions: $a_{i,j,m}=0$ if $i=i_0$ or $j=j_0$,
 and $a_{i_0,j_0,m} = -\frac{1}{\mu_{Z,i_0,j_0}(t_1,t)}$ if $\mu_{Z,i_0,j_0}(t_1,t) \ne 0$, otherwise $a_{i_0,j_0,m} = 0$ (by (3.3) $\mu_{Z,i_0,j_0}(t_1,t) \ne 0$ for all enough large $m$). Then it is not difficult to verify that $\det [I +A_m\mu_Z(t_1,t_m)] = 0$ for all enough large $m$. Hence, by (2.2) $Z_m(t)$ are not $t_1$-regular for all enough large $m$. Moreover, it follows from (3.3) that $\lim\limits_{m \to +\infty} ||Z_m(t_1) - Z(t_1)|| = 0$. Therefore $Z(t)$ is not $t_1$-normal.

Sufficiency. Assume $||\mu_Z(t_1,t)|| \le M, \ph t \ge t_1$ for some $M = const > 0$. Then for every solution $Z_1(t)$ of Eq. (1.1) with $||Z_1(t_1) - Z(t_1)|| < \frac{1}{M}$ the following inequality is valid
$$
\det [I + (Z_1(t_1)  - Z(t_1))\mu_Z(t_1,t)] \ne 0, \phh t \ge t_1.
$$
By virtue of (2.2) from here it follows that $Z_1(t)$ is $t_1$-regular. Hence, $Z(t)$ is $t_1$-normal. The theorem is proved.

In virtue of (2.6) from Theorem 3.1 we obtain immediately.

{\bf Corollary 3.1.} {\it The following statements are valid.

\noindent
1) $t_1$-regular solutions $Z_1(t)$ and $Z_2(t)$ of Eq. (1.1) are $t_1$-normal if and only if the function
$$
I_{Z_1,Z_2}(t) \equiv \il{t_1}{t} Re \hskip 3pt tr  [P(\tau)(Z_1(\tau) - Z_2(\tau))] d\tau, \phh t \ge t_1
$$
is bounded.

\noindent
2) If $Z_N(t)$ and $Z_*(t)$ are respectively  $t_1$-normal and $t_1$-extremal solutions of Eq. (1.1), then
$$
\limsup_{t \to +\infty}\il{t_1}{t} Re \hskip 3pt tr  [P(\tau)(Z_*(\tau) - Z_N(\tau))] d\tau < +\infty,
$$
$$
\liminf_{t \to +\infty}\il{t_1}{t} Re \hskip 3pt tr  [P(\tau)(Z_*(\tau) - Z_N(\tau))] d\tau = -\infty.
$$

\noindent
3) If  $Z_*(t)$ and $Z^*(t)$ are $t_1$-extremal solutions of Eq. (1.1), then
$$
\limsup_{t \to +\infty}\il{t_1}{t} Re \hskip 3pt tr  [P(\tau)(Z_*(\tau) - Z^*(\tau))] d\tau = +\infty,
$$
$$
\liminf_{t \to +\infty}\il{t_1}{t} Re \hskip 3pt tr  [P(\tau)(Z_*(\tau) - Z^*(\tau))] d\tau = -\infty.
$$
}
\phantom{aaaaaaaaaaaaaaaaaaaaaaaaaaaaaaaaaaaaaaaaaaaaaaaaaaaaaaaaaaaaaa} $\blacksquare$

{\bf Definition 3.2.} {\it A regular Eq. (1.1) is called normal, if it has no  extremal solutions.}

{\bf Definition 3.4.} {\it A regular Eq. (1.1) is called irreconcilable, if its every solution is extremal.}

{\bf Definition 3.5.} {\it A regular Eq. (1.1) is called sub extremal, if it has only one extremal solution.}

{\bf Definition 3.5.} {\it A regular Eq. (1.1) is called suer extremal, if it has at least two extremal and other normal solutions.}

From Definitions 3.3 - 3.6 is seen that every regular Eq. (1.1) is or else normal or else irreconcilable or else sub extremal or else super extremal.

For any $t_1$-regular solution $Z(t)$ of Eq. (1.1) denote by $\Omega_Z(t_1)$ the set of all matrices $\Lambda$ of dimension $n\times n$, satisfying the conditions

\noindent
$\alpha$) $\det [I + \Lambda\mu_Z(t_1,t)] \ne 0, \phh t \ge t_1,$

\noindent
$\beta$) there exists a matrix $K$ of dimension $n\times n$, such that $K$ is a limiting matrix for $\mu_Z(t_1,t)$ as $t \to +\infty$ and $\det [I + \Lambda K] = 0.$

{\bf Theorem 3.2.} {\it Let $Z(t)$ be a $t_1$-normal solution of Eq. (1.1). In order that Eq. (1.1) has a $t_1$-extremal solution it is necessary and sufficient that $\Omega_Z(t_1) \ne \emptyset$. If this condition holds, then all $t_1$-extremal solutions $Z_\Lambda(t)$ of Eq. (1.1) are given by the formula
$$
Z_\Lambda(t) = Z(t) + \Psi_Z^{-1}(t)[I + \Lambda \mu_Z(t_1,t)]^{-1}\Lambda\Phi_Z^{-1}(t), \phh t \ge t_1, \ph \Lambda \in \Omega_Z(t_1). \eqno (3.4)
$$
}

Proof. Necessity. Let $Z_*(t)$ be a $t_1$-extremal solution of Eq. (1.1). Then by (2.2) for $\Lambda_* \equiv Z_*(t_1) - Z(t_1)$ the condition $\alpha$) is satisfied. Let $\omega_Z(t_1)$ be the set of all limiting matrices for $\mu_Z(t_1,t)$ as $t \to +\infty$. Since $Z(t)$ is $t_1$-normal  by Theorem 3.1 $\mu_Z(t_1,t)$ is bounded by $t$ on $[t_1,+\infty)$. Then as far as $\mu_Z(t_1,t)$ is continuous $\omega_Z(t_1)$ is a non empty compact set. Show that there exists $K\in \omega_Z(t_1)$ such that
$$
\det [I + \Lambda_* K] = 0.
$$
Suppose this is not so. Then
$$
\delta \equiv \min\limits_{K\in\omega_Z(t_1)}|\det [I + \Lambda K]| > 0.
$$
Let $\cup_{j=1}^N B_j$ be a finite coating with open balls $B_j \ph (j=\overline{1,N})$ of so small radius such that for every $j=\overline{1,N}$ and $K \in B_j$ the inequality
$$
|\det [I + \Lambda_* K]| > \delta/2 \eqno (3.5)
$$
is valid. Note that there exists $T > t_1$ such that
$$
\mu_Z(t_1,t) \subset \cup_{j=1}^N B_j, \phh t \ge T.
$$
(indeed, otherwise there exists a limiting matrix $K$ for $\mu_Z(t_1,t)$ as $t \to +\infty$, not belonging to $\omega_Z(t_1)$). Therefore by (3.5)
$$
|\det [I + \Lambda_*\mu_Z(t_1,t)]| \ge \delta/2, \phh t \ge T. \eqno (3.6)
$$
Let $\delta_1 \equiv \min\limits_{t\in[t_1,T]}|\det [I + \Lambda_*\mu_Z(t_1,t)]| \ge 0$. From here and from (3.6) it follows
$$
|\det [I + \Lambda_* K]| \ge \min\{\delta/2, \delta_1\} > 0, \phh t \ge t_1.
$$
By virtue of Cramer's theorem (the formula of Cramer) from here we obtain that the matrix function $[I + \Lambda_*\mu_Z(t_1,t)]^{-1}, \ph t \ge t_1$ is bounded by $t$ on $[t_1,+\infty)$. Then
$$
r_0 \equiv \bigl(\sup\limits_{t \ge t_1}||[I + \Lambda_*\mu_Z(t_1,t)]^{-1}|| \cdot||\mu_Z(t_1,t)||\bigr)^{-1} > 0.
$$
It follows from here that for every matrix $V$ of dimension $n\times n$ with $||v|| \le r_0/2$ the matrix $I + [I+ \Lambda_*\mu_Z(t_1,t)]^{-1}V\mu_Z(t_1,t)$ is invertible  for all $t\ge t_1$. This together with the easily verifiable equality
$$
I + (\Lambda_* + V)\mu_Z(t_1,t) = [I + \Lambda_* \mu_Z(t_1,t)][I + (I + \Lambda_*\mu_Z(t_1,t))^{-1}V\mu_Z(t_1,t)], \phh t \ge t_1
$$
implies that for every matrix $V$ of dimension $n\times n$ with $||v|| < r_0/2$ the matrix $I + (\Lambda_* + V)\mu_Z(t_1,t)$ is invertible for all $t \ge t_1.$ By (2.2) from here it follows that every solution $\widetilde{Z}(t)$ of Eq. (1.1) with $||\widetilde{Z}(t_1) - Z_*(t_1)|| < r_0/2$ is $t_1$-regular. Hence, $Z_*(t)$ is $t_1$-normal, which contradicts our supposition that $Z_*(t)$ is $t_1$-extremal. The obtained contradiction shows that $\Omega_Z(t_1) \ne \emptyset$.

Sufficiency. Let $\Lambda$ be a matrix of dimension $n\times n$ for which the conditions $\alpha$) and $\beta$) are satisfied, and let $Z_\Lambda(t)$ be the solution of Eq. (1.1) with $Z_\lambda = Z(t_1) + \Lambda$. Then by (2.2)
$$
Z_\Lambda(t) = Z(t) + \Psi_Z^{-1}(t)[I + \Lambda\mu_Z(t_1,t)]^{-1} \Lambda\Phi_Z^{-1}(t), \phh t \ge t_1,
$$
i. e., $Z_\Lambda(t)$ is $t_1$-regular, for which (3.4) is valid. To complete the proof of the theorem it remains to show that $Z_\Lambda(t)$ is $t_1$-extremal. By (2.2) we have
$$
\det \{[I + \Lambda\mu_Z(t_1,t)][I - \Lambda \mu_{Z_\Lambda}(t_1,t)]\} \equiv 1, \phh t \ge t_1. \eqno (3.7)
$$
From the condition $\beta$) it follows that there exists an infinitely large sequence $t_1 < \dots < t_m , \dots$ such that
$$
\lim\limits_{m \to +\infty} \det [I + \Lambda \mu_Z(t_1,t_m)] = 0.
$$
This together with (3.7) implies that
$$
\lim\limits_{m \to +\infty} |\det [I - \Lambda \mu_{Z_\Lambda}(t_1,t_m)]| = +\infty.
$$
Therefore $\mu_{Z_\Lambda}(t_1,t)$ is bounded by $t$ on $[t_1,+\infty)$. In virtue of Theorem 3.1 from here it follows that $Z_\Lambda(t)$ is $t_1$-extremal. The theorem is proved.

For any $t_1$-regular solution $Z(t)$ of Eq. (1.1.) set
$$
\nu_Z(t)\equiv \il{t}{+\infty}\Phi_Z^{-1}(\tau) P(\tau) \Psi_Z^{-1}(\tau) d \tau,\phh t \ge t_1,
$$
where $ \Phi_Z(t)$ and $\Psi_Z(t)$ are solutions of the linear matrix equations
$$
\Phi' = [P(t)Z(t) + R(t)]\Phi = 0, \phh t \ge t_1,
$$
$$
\Psi' = \Psi[Z(t) P(t) + Q(t)] = 0, \phh t \ge t_1
$$
respectively with $\Phi_Z(t_1) = \Psi_Z(t_1) = I.$ Note that if $\nu_Z(t)$ converges for some $t=t_2\ge t_1$, then it converges for all $t \ge t_1$.
\pagebreak

{\bf Theorem 3.3.} {\it Let for some $t_1$-regular solution $Z_0(t)$ of Eq. (1.1) the integral $\nu_{Z_0}(t_1)$ be convergent and $\det \nu_{Z_0}(t) \ne 0, \ph t \ge t_1$. Then the following statements are valid.

\noindent
a) The matrix function
$$
Z_*(t) \equiv Z_0(t) - [\Phi_{Z_0}(t) \nu_{Z_0}(t)\Psi_{Z_0}(t)]^{-1}, \phh t \ge t_1
$$
is a $t_1$-extremal solution of Eq. (1.1).

\noindent
b)
$$
\nu_{Z_*}(t) = \infty, \phh t \ge t_1.
$$

\noindent
c) For every $t_1$-normal solution $Z_N(t)$ of Eq. (1.1) the equality
$$
\il{t_1}{+\infty} Re \hskip 2pt tr \hskip 2pt [P(\tau)(Z_*(\tau) - Z_N(\tau))]d \tau = -\infty \eqno (3.8)
$$
is valid.}

Proof. Let $Z_*(t)$ be the solution of Eq. (1.1) with $Z_*(t_1) = Z_0(t_1) -\linebreak -[\Phi_{Z_0}(t_1)\nu_{Z_0}(t_1)\Psi_{Z_0}(t_1)]^{-1} = Z_0(t_1) - [\nu_{Z_0}(t_1)]^{-1}$, and let $\Lambda\equiv Z_*(t_1) - Z_0(t_1) =  \linebreak -[\nu_{Z_0}(t_1)]^{-1}.$ Show that
$$
\det [I + \Lambda \mu_{Z_0}(t_1,t)] \ne 0, \phh t \ge t_1.
$$
Suppose for some $t_2 \ge t_1$
$$
\det [I + \Lambda \mu_{Z_0}(t_1,t_2)] \ne 0. \eqno (3.9)
$$
Then taking into account the equality $\mu_{Z_0}(t_1,t_2) = \nu_{Z_0}(t_1)- \nu_{Z_0}(t_2)$ we obtain
$$
0 = \det [I + \Lambda \mu_{Z_0}(t_1,t_2)] = \det [I - \nu_{Z_0}(t_1)]^{-1} (\nu_{Z_0}(t_1) - \nu_{Z_0}(t_2)) = \det \nu_{Z_0}(t_2) \ne 0.
$$
We have obtained a contradiction, which proves (3.9). By virtue of (2.2) from (3.9) it follows that $Z_*(t)$ is $t_1$-regular. Since $\lim\limits_{t \to +\infty}[I + \Lambda \mu_{Z_0}(t_1,t)] = \lim\limits_{t \to +\infty}\nu_{Z_0}(t) = 0$ by Theorem~ 3.2 $Z_*(t)$ is $t_1$-extremal and
$$
Z_*(t) = Z_0(t) + \Psi_{Z_0}^{-1}(t)[I - \nu_{Z_0}^{-1}(t_1)\mu_{Z_0}(t_1,t)]^{-1} [\nu_{Z_0}(t_1)]^{-1} \Phi_{Z_0}^{-1}(t)  = \phantom{aaaaaaaaaaaaaaaaaaa}
$$
$$
= Z_0(t) + \Psi_{Z_0}^{-1}(t)[I - \nu_{Z_0}^{-1}(t_1)(\nu_{Z_0}(t_1) - \nu_{Z_0}(t))]^{-1} [\nu_{Z_0}(t_1)]^{-1} \Phi_{Z_0}^{-1}(t) =
$$
$$
\phantom{aaaaaaaaaaaaaaaaaaaaaaaaaaaaaaaaaaaaa}= Z_0(t)  - [\Phi_{Z_0}(t) \nu_{Z_0}(t)\Psi_{Z_0}(t)]^{-1}, \phh t \ge t_1.
$$
The statement a) is proved. Prove b). By (2.7) we have
$$
\det \{[I + \nu_{Z_0}^{-1}(t_1)\mu_{Z_*}(t_1,t)] \hskip 1pt [I - \nu_{Z_0}^{-1}(t_1) \mu_{Z_0}(t_1,t)]\} \equiv 1, \phh t \ge t_1.
$$
From here and from the relation
$$
\lim\limits_{t\to +\infty}\det [I - \nu_{Z_0}^{-1}(t_1) \mu_{Z_0}(t_1,t)] = 0 \eqno (3.10)
$$
it follows that
$$
\lim\limits_{t\to +\infty}|\det \{[I + \nu_{Z_0}^{-1}(t_1)\mu_{Z_*}(t_1,t)]| = +\infty.
$$
Hence, $\nu_{Z_*}(t_1) = \infty$. Therefore, $\nu_{Z_*}(t) = \infty$ for all $t \ge t_1$. The statement b) is proved. Finally by (2.6) from (3.10) it follows (3.8). The theorem is proved.
\vskip 20 pt

{\bf 4. The asymptotic behavior of solutions of linear systems of ordinary differential equations.} Let $A(t) \ph B(t), \ph C(t)$ and $D(t)$ be complex-valued matrix functions of dimension $n\times n$ on $[t_0,+\infty)$. Consider the linar systems of matrix ordinary differential equations
$$
\sist{\Phi' = A(t)\Phi + B(t)\Psi,}{\Psi' = C(t)\Phi + D(t)\Psi, \ph t \ge t_0.} \eqno (4.1)
$$
Here $\Phi = \Phi(t)$ and $\Psi = \Psi(t)$ are unknown continuously differentiable  matrix functions of dimension $n\times n$ on $[t_0,+\infty)$. In (4.1) substitute
$$
\Psi = Z\Phi. \eqno (4.2)
$$
We obtain
$$
\sist{\Phi' = [B(t) Z + A(t)]\Phi,}{[Z' + Z B(t) Z + Z A(t) - D(t)Z - C(t)]\Phi = 0, \ph t \ge t_0.}
$$
It follows from here that all solutions $Z(t)$ of the matrix Riccati equation
$$
Z' + Z B(t) Z + Z A(t) - D(t) Z - C(t) = 0, \phh t \ge t_0, \eqno (4.3)
$$
existing on $[t_1,t_2) \ph (t_0 \le t_1 < t_2 \le +\infty)$ are connected with solutions $(\Phi(t),\Psi(t))$ of the system (4.1) by the relations
$$
\Phi'(t) = [B(t) Z(t) + A(t)] \Phi(t), \ph \Phi(t_1) \ne 0, \phh \Psi(t) = Z(t) \Phi(t), \phh t\in [t_1,t_2). \eqno (4.4)
$$

{\bf Definition 4.1} {\it A solution $(\Phi(t),\Psi(t))$ of the system (4.1) is called $t_1$-regular $(t_1 \ge t_0)$ if $\det \Phi(t) \ne 0, \ph t \ge t_1.$}

{\bf Definition 4.2.} {\it A $t_1$-regular solution $(\Phi(t),\Psi(t))$ of the system (4.1) is called principal (non principal) if $Z(t) \equiv \Psi(t)\Phi^{-1}(t), \ph t \ge t_1$ is a $t_1$-extremal ($t_1$-normal) solution of Eq.~ (4.3).}

{\bf Definition 4.3.} {\it The system (4.1) is called regular, if it has at least one regular solution.}

{\bf Definition 4.4.} {\it The regular system (4.1) is called normal (irreconcilable, sub extremal, super extremal, extremal) if Eq. (4.3) is normal (irreconcilable, sub extremal, super \linebreak extremal, extremal).}

Let $(\Phi(t),\Psi(t))$ be a $t_1$-regular solution of the system (4.1) and let $Z(t)\equiv \Psi(t)\Phi^{-1}(t), \linebreak t \ge t_1$. Then by (4.4) and the Liouville's formula we have
$$
|\det \Phi(t)| = |\det \Phi(t_1)|\exp\biggl\{\il{t_1}{t}Re \hskip 2pt tr[B(\tau) Z(\tau) + A(\tau)] d \tau\biggr\}, \phh t \ge~ t_1.
$$
By Corollary 3.1 from here we obtain immediately

{\bf Theorem 4.1.} {\it The following statements are valid.

\noindent
I) If the system (4.1) is normal, then for its two arbitrary regular solutions $(\Phi_m(t),\Psi_m(t)), \linebreak m=1,2$ the relations
$$
\limsup\limits_{t \to +\infty}\frac{|\det\Phi_1(t)|}{|\det\Phi_2(t)|} < +\infty, \phh \limsup\limits_{t \to +\infty}\frac{|\det\Phi_2(t)|}{|\det\Phi_1(t)|} < +\infty
$$
are fulfilled.

\noindent
II) If the system (4.1) is irreconcilable, then for is two arbitrary linearly independent regular solutions $(\Phi_m(t),\Psi_m(t)), \ph m=1,2$ the equalities
$$
\limsup\limits_{t \to +\infty}\frac{|\det\Phi_1(t)|}{|\det\Phi_2(t)|} = \limsup\limits_{t \to +\infty}\frac{|\det\Phi_2(t)|}{|\det\Phi_1(t)|} = +\infty
$$
are fulfilled.

\noindent
III) If the system (4.1) is sub extremal, then there exists a regular solution $(\Phi_*(t),\Psi_*(t))$ of  (4.1) such that, for every regular solutions $(\Phi_m(t),\Psi_m(t)), \ph m=1,2$, linearly independent of $(\Phi_*(t),\Psi_*(t))$, the relations
$$
\limsup\limits_{t \to +\infty}\frac{|\det\Phi_*(t)|}{|\det\Phi_1(t)|} < +\infty, \phh \liminf\limits_{t \to +\infty}\frac{|\det\Phi_*(t)|}{|\det\Phi_1(t)|} = 0,
$$
$$
\limsup\limits_{t \to +\infty}\frac{|\det\Phi_1(t)|}{|\det\Phi_2(t)|} < +\infty, \phh \limsup\limits_{t \to +\infty}\frac{|\det\Phi_2(t)|}{|\det\Phi_1(t)|} < +\infty
$$
are fulfilled.

\noindent
IV) If the system (4.1) is super extremal, then there exists its two regular solutions $(\Phi_*(t),\Psi_*(t))$ and $(\Phi^*(t),\Psi^*(t))$ such that
$$
\limsup\limits_{t \to +\infty}\frac{|\det\Phi_*(t)|}{|\det\Phi^*(t)|} = \limsup\limits_{t \to +\infty}\frac{|\det\Phi^*(t)|}{|\det\Phi_*(t)|} = +\infty
$$
and for all its two arbitrary solutions  $(\Phi_m(t),\Psi_m(t)), \ph m=1,2$, linearly independent of $(\Phi_*(t),\Psi_*(t))$ and $(\Phi^*(t),\Psi^*(t))$ the following relations are fulfilled:
$$
\limsup\limits_{t \to +\infty}\frac{|\det\Phi_1(t)|}{|\det\Phi_2(t)|} < +\infty, \phh \limsup\limits_{t \to +\infty}\frac{|\det\Phi_2(t)|}{|\det\Phi_1(t)|} < +\infty,
$$
$$
\limsup\limits_{t \to +\infty}\frac{|\det\Phi_*(t)|}{|\det\Phi_m(t)|} < +\infty, \phh \limsup\limits_{t \to +\infty}\frac{|\det\Phi^*(t)|}{|\det\Phi_m(t)|} < +\infty,
$$
$$
\liminf\limits_{t \to +\infty}\frac{|\det\Phi_*(t)|}{|\det\Phi_m(t)|} = \liminf\limits_{t \to +\infty}\frac{|\det\Phi^*(t)|}{|\det\Phi_m(t)|} = +\infty, \ph m=1,2.
$$
}
\phantom{aaaaaaaaaaaaaaaaaaaaaaaaaaaaaaaaaaaaaaaaaaaaaaaaaaaaaaaaaaaaaaaaaaa} $\blacksquare$

\vskip 20pt

\centerline{\bf References}

\vskip 20pt

\noindent
1. G. Freiling, A survey  of nonsymmetric Riccati equations. Linear Algebra and its Applications.  351-352 (2002) pp. 243-270.

\noindent
2. H. W. Knobloch, M.Pohl, On Riccati matrix differential equations. Results Math. 31 (19970, PP. 337-364.

\noindent
3. G. Freiling, G. Jank, A. Sarychev, Non-blow-up conditions for Riccati-type matrix differential and difference equations. Results math. 37 (1998) pp. 84-103.

\noindent
4. G. A. Grigorian, On some properties of solutions of the Riccati equation.
Izvestiya  NAS \linebreak\phantom{aa} of Armenia, vol. 42, $N^\circ$ 4, 2007, pp. 11 - 26.

\noindent
5. G. A. Grigorian,    Criteria of global solvability for Riccati scalar equations.
Izv. Vyssh.  \linebreak\phantom{aa} Uchebn. Zaved. Mat., 2015, Number 3, pp. 35–48.

\noindent
6.  G. A. Grigorian, Some properties of differential root and their applications.  Acta Math.  \linebreak\phantom{aa} Univ. Comenianae
Vol. LXXXV,2 (2016) pp. 205-217.

\noindent
7. G. A. Grigorian, Properties of solutions of the scalar Riccati equations with complex  \linebreak\phantom{aa} coefficients and some their applications, Diff Eq. and Appl, vol 10, Num. 3 (2018) pp. 277-298.

\noindent
8. G. A. Grigorian, Global solvability criteria for quaternionic Riccati equations. Archivum \linebreak\phantom{aa}   Mathematicum. Tomus 57 (2021), pp. 83–99.

\noindent
9.  G. A. Grigorian, Properties of solutions of quaternionic
Riccati equations, Archivum  \linebreak\phantom{aa}  Mathematicum. In print.

\noindent
10. G. A. Grigorian, On the Stability of Systems of Two First-Order Linear Ordinary  \linebreak\phantom{aa}
Differential Equations, Differ. Uravn., 2015, vol. 51, no. 3, pp. 283 - 292.

\noindent
11. G. A. Grigorian,   Stability Criterion for Systems of Two First-Order Linear Ordinary   \linebreak\phantom{aa}   Differential Equations.
Mat. Zametki, 2018, Volume 103, Issue 6, Pp. 831–840.

\noindent
12.  G. A. Grigorian, On one oscillatory criterion for the second order linear
 ordinary \linebreak \phantom{a} differential equations. Opuscula Math. 36, no. 5 (2016), 589–601. \linebreak \phantom{aa}
   http://dx.doi.org/10.7494/OpMath.2016.36.5.589

\noindent
13. G. A. Grigorian, Some properties of the solutions of third order linear ordinary \linebreak \phantom{a}    differential  equations.  Rocky Mountain Journal of Mathematics, vol. 46, no. 1,  2016, \linebreak \phantom{a}  pp. 147 - 161.

\noindent
14. G. A. Grigorian,    Oscillatory criteria for the second order linear ordinary differential \linebreak \phantom{a} equations.
Math. Slovaca 69 (2019), No. 4, pp. 857-870.

\noindent
15. G. A. Grigorian,  Oscillatory criteria for the systems of two  first order
 linear ordinary  \linebreak\phantom{aa} differential equations. Rocky J. Math. Vol. 47, Number 5, 2017, pp. 1497-1534.

\noindent
16. G. A. Grigorian,  Oscillatory and non oscillatory criteria for the systems of two linear  \linebreak\phantom{aa}
first order two by two dimensional matrix ordinary differential equations. Arch Math. \linebreak\phantom{aa}
Tomus 54 (2018), pp. 189–203

\noindent
17. A. I. Egorov. Riccati equations. Moskow, Fizmatlit,  2001.

\noindent
18. G. A. Grigorian, A new oscillation criterion for the generalized Hill's equation, Diff  \linebreak\phantom{aa}  Eq. and Appl, vol. 9, Num. 3 (2011) pp. 369-377.

\end{document}